\documentclass[12 pt]{amsart}

\usepackage{graphicx}
\usepackage{color}
\usepackage{amssymb, amsmath, amsthm}
\usepackage{mathptmx}

\usepackage{epsfig}
\usepackage{epstopdf}
\usepackage{graphicx}

\newtheorem{thm}{Theorem}[section]

\newtheorem{cor}[thm]{Corollary}

\newtheorem{defn}[thm]{Definition}
\newtheorem{lemma}[thm]{Lemma}

\newtheorem{conj}{Conjecture}

\newcommand{\R}{\mathbb R} %REALS

\newcommand{\bi}{\begin{itemize}}
\newcommand{\ei}{\end{itemize}}
\newcommand{\be}{\begin{enumerate}}
\newcommand{\ee}{\end{enumerate}}

\newcommand{\n}{\beta}

\newcommand{\emp}{\emptyset}
\newcommand{\X}{\times}

\newcommand{\eps}{\epsilon}

\newcommand{\M}{\mathcal{M}}

\newcommand{\A}{\alpha}
\newcommand{\pd}{\partial}

\newcommand{\MP}{\mathcal{MP}}

\newcommand{\K}{\mathcal{K}}

\newcommand{\Ss}{\mathcal{S}}
\newcommand{\BT}{\mathcal{BT}}
\newcommand{\KT}{\mathcal{KT}}
\newcommand{\LL}{\mathcal{L}}

\begin{document}
\title{Properties of knots preserved by cabling}
\author{Alexander Zupan}

\maketitle
\begin{abstract}
We examine geometric properties of a knot $J$ that are unchanged by taking a $(p,q)$-cable $K$ of $J$.  Specifically, we show that $w(K) = q^2 \cdot w(J)$, where $w(K)$ is the width of $K$ in the sense of Gabai.  We use this fact to demonstrate that thin position is a minimal bridge position of $J$ if and only if the same is true for $K$, and more generally we show that any thin position of $K$ is an ``obvious" cabling of a thin position of $J$.  We conclude by proving that $J$ is meridionally small (mp-small) if and only if $K$ is meridionally small (mp-small), and if $J$ is mp-small and every non-minimal bridge position for $J$ is stabilized, then the same is true for $K$.
\end{abstract}

\section{Introduction}
In \cite{schub}, Horst Schubert introduces the bridge number $b(K)$ of a knot $K$ in $S^3$ and proves that for a satellite knot $K$ with companion $J$ and pattern $\hat{K}$ with wrapping number $n$, $b(K) \geq n \cdot b(J)$.  This is also proven by Schultens in \cite{bridge}.  In particular, if $K$ is a $(p,q)$-cable of $J$, then Schubert's bound is tight; that is, $b(K) = q \cdot b(J)$.  Bridge number can be approached as an element of the broad collection of topological invariants that can be calculated by minimizing some sort of numerical complexity over all possible positions of a space with respect to a Morse function.  Knot width, introduced by Gabai in \cite{gabai}, also falls into this category, as do Heegaard genus and thin position of 3-manifolds.  In fact, knot width can be seen as a type of loose refinement of bridge number: for any knot $K$ whose exterior contains no essential meridional planar surfaces, $w(K) = 2 \,b(K)^2$ and any thin position is a minimal bridge position for $K$.  We call any knot satisfying $w(K) = 2 \, b(K)^2$ \emph{bridge-thin} and denote the collection of all such knots $\BT$.  For knots $K$ not in $\BT$, the situation becomes more complicated: recently Blair and Tomova have shown that the bridge number of a knot $K$ cannot always be recovered from the thin position of $K$ \cite{blairtom}. \\

In order to further develop the connection between bridge number and knot width, we make the following conjecture in \cite{zupan1} as an analogue to the above theorem of Schubert and Schultens:
\begin{conj}
If $K$ is a satellite knot with companion $J$ and pattern $\hat{K}$ with wrapping number $n$, then
\[ w(K) \geq n^2 \cdot w(J).\]
\end{conj}
One purpose of this paper is to prove the conjecture for cable knots and to show that the bound is tight, as with bridge number.  We prove
\begin{thm}
If $K$ is a $(p,q)$-cable of $J$, then
\[ w(K) = q^2 \cdot w(J),\]
and any thin position of $K$ is an ``obvious" cable of a thin position of $J$.
\end{thm}
In particular, this shows that $J$ is in $\BT$ if and only if some cable of $J$ is in $\BT$. \\

The collection $\BT$ contains many knots.  Let $\MP$ denote the collection of knots $K$ whose exterior $E(K)$ contains no essential meridional planar surface (mp-small knots), let $\M$ denote the collection whose exterior $E(K)$ contains no essential meridional surface (meridionally small knots), and let $\Ss$ consist of those knots whose exterior $E(K)$ contains no essential closed surface (small knots).  By \cite{dehn}, $\Ss \subset \M$, clearly $\M \subset \MP$, and by \cite{thomps}, $\MP \subset \BT$.  In summary,
\[ \Ss \subset \M \subset \MP \subset \BT.\]
It is well known that $\Ss \neq \M$; for instance, see Proposition 1.6 of \cite{morim}.  Examples of knots $K \in \BT$ such that $K \notin \MP$ appear in \cite{blair} and \cite{heath}.  We are not aware of examples of knots $K$ in $\MP$ and not in $\M$, but we suspect that such examples exist. \\

For any knot $J$ and $(p,q)$-cable $K$ of $J$, we have that $K \notin \Ss$ as the companion torus is essential in $E(K)$, so smallness is not preserved by cabling.  However, since cabling does preserve containment in $\BT$, we can examine the situation when $J \in \M$ or $J \in \MP$.  We prove
\begin{thm}
For a knot $J$, the following are equivalent:
\be
\item[(a)] $J \in \M$.
\item[(b)] Every cable $K$ of $J$ is in $\M$.
\item[(c)] There exists a cable $K$ of $J$ in $\M$.
\ee
The same statement is true if $\MP$ replaces $\M$.
\end{thm}
Thus, meridional smallness and mp-smallness are also preserved by cabling. \\

Finally, a natural problem concerning bridge number is to determine for which knots are all non-minimal bridge positions stabilized.  A bridge position is stabilized if a minimum and maximum can be cancelled to create a presentation with a lower bridge number.  In this direction, Otal shows in \cite{otal} that every non-minimal bridge position of the unknot or a 2-bridge knot is stabilized, and Ozawa has recently proved the same statement for torus knots in \cite{ozawa}.  We follow Ozawa's proof in order to show the following:
\begin{thm}
Suppose $K$ is a $(p,q)$-cable of $J$, where $J$ is mp-small.
\be
\item[(a)] If every non-minimal bridge position of $J$ is stabilized, then every non-minimal bridge position of $K$ is stabilized.
\item[(b)] The cardinality of the collection of minimal bridge positions of $K$ does not exceed the cardinality of the collection of minimal bridge positions of $J$.
\ee
\end{thm}

\section{Preliminaries}
Fix a Morse function $h:S^3 \rightarrow \R$ such that $h$ has exactly two critical points, which we denote $\pm \infty$.  Now, fix a knot $K$ and let $\mathcal{K}$ denote the set of all embeddings $k$ of $S^1$ into $S^3$ isotopic to $K$ and such that $h \mid_k$ is Morse.  For each such $k$, there are critical values $c_0 < \dots < c_p$ of $h \mid_k$.  Choose regular values $c_0 < r_1 < \dots < r_p < c_p$ of $h$, and for each $i$, let $x_i = |k \cap h^{-1}(r_i)|$; thus, we associate the tuple $(x_1,\dots,x_p)$ of even integers to $k$.  Define the width of $k$ to be
\[ w(k) = \sum x_i\]
and the bridge number of $k$, $b(k)$, to be the number of maxima (or minima) of $h\mid_k$.  Then the width and bridge number of the knot $K$ are defined as
\[ w(K) = \min_{k \in \mathcal{K}} w(k) \qquad \text{ and } \qquad b(K) = \min_{k \in \mathcal{K}} b(k).\]
Any $k \in \mathcal{K}$ satisfying $w(k) = w(K)$ is a called a \emph{thin position} of $K$, and any $k$ with all maxima above all minima is called a \emph{bridge position} of $K$.  If $k$ is a bridge position and $b(k) = b(K)$, we call $k$ a \emph{minimal bridge position} of $K$. \\

For any $k \in \mathcal{K}$ with associated tuple $(x_1,\dots,x_p)$, we say $x_i$ corresponds to a thick level if $x_i > x_{i-1},x_{i+1}$, and $x_i$ corresponds to a thin level if $x_i < x_{i-1},x_{i+1}$, where $2 \leq i \leq n-1$.  Thus, we can associate a thick/thin tuple of integers $(a_1,b_1,a_2,\dots,b_{n-1},a_n)$ to $k$, where each $a_i$ is corresponds to a thick level and each $b_i$ corresponds to a thin level.  We also have a collection of level surfaces $\hat{A}_i$ and $\hat{B}_i$ satisfying $|k \cap \hat{A}_i| = a_i$ and $|k \cap \hat{B}_i| = b_i$.  We call these surfaces thick and thin surfaces, respectively, for $k$.  From \cite{planar},
\begin{equation}\label{wsquare}
w(k) = \frac{1}{2} \left( \sum a_i^2 - \sum b_i^2 \right).
\end{equation}
The definitions in the next two paragraphs are taken from \cite{width}.  Given an embedding $k$, there may exist certain isotopies that decrease $w(k)$, as determined by the intersections of strict upper and lower disks.  An \emph{upper disk} for $k$ at a thick surface $\hat{A}_i$ is an embedded disk $D$ such that $\pd D$ consists of two arcs, one arc in $k$ and one arc in $\hat{A}_i$, where the arc in $k$ does not intersect any thin surfaces and contains exactly one maximum.  A \emph{strict upper disk} is an upper disk whose interior contains no critical points with respect to $h$.  A \emph{lower disk} and \emph{strict lower disk} are defined similarly. \\

If there is a pair $(D,E)$ of strict upper and lower disks for $k$ at $\hat{A}_i$ that intersect in a single point contained in $k$, we can cancel out the maximum and minimum of $k$ contained in $D$ and $E$, and $w(k)$ decreases by $2a_i - 2$.  We call this a \emph{type I move}.  If $k$ admits a type I move at $\hat{A}_i$, we say $k$ is \emph{stabilized} at $\hat{A}_i$.  If $D$ and $E$ are disjoint, we can slide the minimum above the maximum with respect to $h$, decreasing $w(k)$ by 4.  This is called a \emph{type II move}. \\

\section{Cable knots and the companion torus}
The knots we will study are called cable knots, a specific type of satellite knot, defined below:
\begin{defn}
Suppose that $\hat{K}$ is a knot contained in a solid torus $V$ with core $C$ such that every meridian of $V$ intersects $\hat{K}$, and let $J$ be any nontrivial knot in $S^3$.  Further, let $\varphi:V \rightarrow S^3$ be an embedding such that $\varphi(C)$ is isotopic to $J$.  Then we say that $K = \varphi(\hat{K})$ is a \textbf{satellite knot} with companion $J$ and pattern $\hat{K}$.  We call $\pd \varphi(V)$ the \textbf{companion torus} corresponding to $J$.
\end{defn}
In general, if $K$ is a satellite knot, $K$ might not be isotopic into the companion torus.  However, if its pattern $\hat{K}$ is a torus knot, we can push the pattern $\hat{K}$ into $T = \pd \varphi(V)$, an assumption that yields more information than the general case when we study the foliation of $T$ induced by $h$.  This leads to the definition of cable knots.
\begin{defn}
A \textbf{$(p,q)$-cable knot} is a satellite knot with pattern a $(p,q)$-torus knot $\hat{K}$.
\end{defn}
\begin{figure}[h!]
  \centering
    \includegraphics[width=0.525\textwidth]{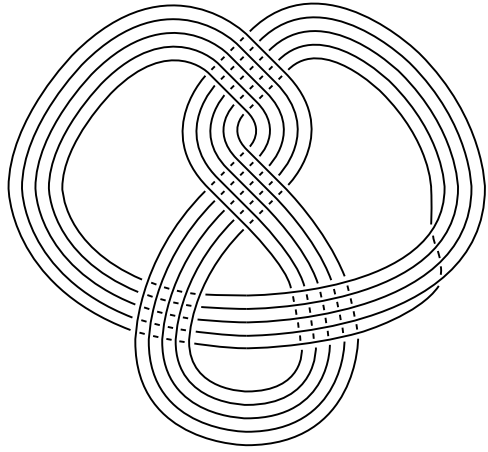}
    \caption{A $(1,5)$-cable of the figure eight knot with the blackboard framing}
\end{figure}
We follow the convention of \cite{lither}, where a $(p,q)$-torus knot in $T$ has homology class $p[\mu]+ q[\lambda]$ with $\mu$ a meridian and $\lambda$ a longitude of $\varphi(V)$.  Observe that a $(p,q)$-torus knot has intersection number $q$ with any meridian of $\varphi(V)$. \\

As noted above, $h$ induces a singular foliation, $F_S$, of any surface $S \subset S^3$ such that $h \mid_S$ is Morse.  As in \cite{bridge}, we distinguish between the various types of saddle points of $F_S$:
\begin{defn}
Let $c$ be a critical value corresponding to a saddle $\sigma$ of $F_S$, so that some component of $S \cap h^{-1}([c-\eps,c+\eps])$ is a pair of pants $P$.  Let $s_1$ and $s_2$ denote the components of $\pd P$ lying in the same level.  If either bounds a disk in $S$, we say that $\sigma$ is an \textbf{inessential saddle}.  If $\sigma$ is not inessential, then it is an \textbf{essential saddle}.
\end{defn}
In addition to type I and type II moves, we will require one more move to systematically simplify cable knots: disk slides.
\begin{defn}
Suppose that some embedding $k$ of a knot $K$ is contained in a surface $S$ such that $h \mid_S$ is Morse.  Suppose further that $D \subset S$ is a disk satisfying the following conditions:
\be
\item $\pd D$ is contained in some level surface $h^{-1}(r)$;
\item $\text{int}(D)$ contains exactly one minimum or maximum;
\item some component $k \cap D$ is an outermost arc $\gamma$ that contains exactly one critical point of $h\mid_k$ and cobounds a disk $\Delta$ with an arc $\pd D$, where $\text{int}(\Delta)$ contains the critical point of $D$.
\ee
Then there is an isotopy through $\Delta$ supported in a neighborhood of $\Delta$ that takes $\gamma$ to an arc $\gamma'$ that cobounds a disk $\Delta'$ with an arc in $\pd D$ such that $\text{int}(\Delta')$ contains no critical points and $\text{int}(\Delta') \cap k = \emp$.  Additionally, $\gamma'$ can be chosen so that it contains exactly one critical point occuring at the same height as the critical point of $\gamma$.  We call the isotopy that replaces $\gamma$ with $\gamma'$ a \textbf{disk slide}, pictured in Figure 2.
\end{defn}
\begin{figure}[h!]\label{disk}
  \centering
    \includegraphics[width=0.675\textwidth]{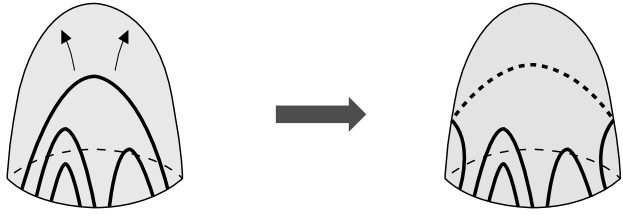}
    \caption{A disk slide}
\end{figure}
Note that by requiring that the critical point of $\gamma'$ occur on the same level as the critical point of $\gamma$, we ensure that if $k'$ is the result of performing a disk slide on $k$, we have $w(k') = w(k)$, although the width of $k$ may temporarily increase during some intermediate step of the disk slide.

\section{Efficient position of $k$}
\label{eff}

From this point forward, we set the convention that $K$ is a $(p,q)$-cable knot with companion $J$ and companion torus $T$.  As noted above, for any $k \in \K$, we may assume there is a torus isotopic to $T$ containing $k$.  In order to study each embedding $k \in \K$ with respect to various representatives in the isotopy class of the companion torus $T$, we define the following collection:
\[ \KT = \left\{ (k,T_k): k \in \K, \, T_k \sim T, \, k \subset T_k, \text{ and } h \mid_{T_k} \text{ is Morse}\right\}.\]
For an arbitrary element $(k,T_k)$ of $\KT$, the torus $T_k$ is compressible only on one side (see Section \ref{small} for a definition), and we let $V_k$ denote the solid torus bounded by $T_k$.  We may also assume that the critical points of $k$ and $T_k$ occur at different levels.  The goal of this section is to show that any embedding $k$ can be deformed to lie on a torus isotopic to $T$ in an efficient way and without increasing $w(k)$.  To this end, consider an arc $\gamma$ of $k$ embedded in an annulus $A$, where both boundary components of $A$ are level and $\text{int}(A)$ contains no critical points with respect to $h$.  Suppose that $\gamma$ is an essential arc containing critical points, and let $x$ denote the lowest maximum of $\gamma$.  There are two points $y,z \in \gamma$ corresponding to minima (one could be an endpoint) such that $\gamma$ contains a monotone arc connecting $x$ to $y$ and $x$ to $z$.  Without loss of generality, suppose that $h(y) > h(z)$.  Then level arc components of $A \cap h^{-1}(h(y) + \eps)$ cobound disks $D,E \subset A$ with arcs in $\gamma$ such that $\pd D$ contains exactly one maximum ($x$), $\pd E$ contains exactly one minimum ($y$), and such that $\text{int}(D) \cap \gamma = \text{int}(E) \cap \gamma =\emp$.  Further, $D \cap E$ is a single point contained in $\gamma$.  Refer to Figure 3. \\
\begin{figure}[h!]\label{arc}
  \centering
    \includegraphics[width=0.225\textwidth]{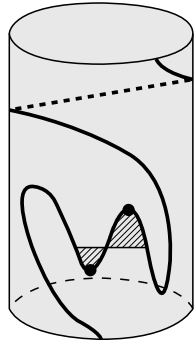}
    \caption{An essential arc $\gamma$ containing critical points, shown with upper and lower disks $D$ and $E$ described above.}
\end{figure}

If there is exactly one thick surface between $x$ and $y$, then $k$ admits a type I move canceling the two critical points.  If there is more than one such thick surface, we can slide $x$ down along $D$, performing some number of type II moves, until there is one such thick surface, after which $k$ admits a type I move.  A similar argument is valid if $\gamma$ is inessential in $A$ or if $\gamma$ is embedded in a disk $\hat{D}$ with level boundary and one critical point, and $\gamma$ contains more than one critical point.  We will use these ideas in the proof of the next lemma, which is similar to the proof of Lemma 3.3 in \cite{zupan1} modified slightly to accommodate the fact that for any $(k,T_k) \in \KT$, we have $k \subset T_k$.
\begin{lemma}
\label{lem2}
For any $k \in \K$, there exists $(k',T_{k'}) \in \KT$ such that $w(k') \leq w(k)$ and the foliation of $T_{k'}$ by $h$ contains no inessential saddles.
\begin{proof}
Suppose that $(k,T_k) \in \KT$ and $F_{T_k}$ contains an inessential saddle.  By Lemma 1 (the Pop Over Lemma) of \cite{bridge}, we can modify $T_k$ so that there exists a saddle $\sigma$ of $F_{T_k}$ with corresponding critical value $c$ such that some component $P$ of $T_k \cap h^{-1}([c-\eps,c+\eps])$ is a pair of pants, with boundary components $s_1,s_2,s_3$ satisfying
{\be
\item $s_1$ and $s_2$ are contained in the same level surface $L$ of $h$,
\item $s_1$ bounds a disk $D_1 \subset T_k$ such that $F_{T_k}$ restricted to $D_1$ contains only one maximum or minimum,
\item $D_1$ co-bounds a 3-ball $B$ with a disk $\tilde{D}_1 \subset L$ such that $B$ does not contain $\pm \infty$ and such that $s_2$ lies outside of $\tilde{D}_1$.
\ee}
Without loss of generality, let $D_1$ contain one maximum.  By the above argument, we can perform type I and II moves so that each component of $k \cap D_1$ contains exactly one maximum.  Let $A = P \cup D_1$, so that $A$ is an annulus with boundary components $s_2$ and $s_3$.  We can choose $P$ so that $k \cap P$ consists only of vertical arcs.  Then if $k \cap D_1 \neq \emp$, $k \cap A$ contains inessential arcs with exactly one critical point and boundary in $s_3$.  Let $\gamma$ denote an outermost arc in $A$, so that $\gamma$ cobounds a disk $\Delta$ with an arc in $s_3$ and $\text{int}(\Delta) \cap k = \emp$.  If $\text{int}(\Delta)$ contains a critical point, we can perform a disk slide to remove it, after which can remove $\gamma$ from $k \cap A$, possibly by type II moves. \\

Thus, after some combination of type I moves, type II moves, and disk slides, we may assume that $k \cap D_1 = \emp$.  By the proof of Lemma 3.3 in \cite{zupan1}, we can cancel the saddle $\sigma$ with the maximum contained in $D_1$ without increasing $w(k)$.  Repeating this process finitely many times finishes the proof.

\end{proof}
\end{lemma}
Hence, for the purpose of finding thin position we may assume that for any $(k,T_k) \in \KT$, the foliation $F_{T_k}$ contains no inessential saddles, and so each disk $D \subset T_k$ with level boundary must contain exactly one critical point.  We distinguish between two types of saddles:
\begin{defn}
Let $(k,T_k) \in \KT$, with $\sigma$ a saddle point of $F_{T_k}$ and $c$ the corresponding critical value.  Then some component $P$ of $T \cap h^{-1}([c - \eps,c + \eps])$ is a pair of pants.  We say that $\sigma$ is a \textbf{lower saddle} if $P \cap h^{-1}(c+\eps)$ has two components; otherwise, $\sigma$ is an \textbf{upper saddle}.
\end{defn}
Next, for $(k,T_k) \in \KT$, we suppose that the foliation of $T_k$ contains only essential saddles and decompose $T_k$ into annuli as follows:  for each critical value $c_i$ corresponding to a saddle, some component of $T_k \cap h^{-1}([c_i - \eps,c_i+\eps])$ is a pair of pants, call it $P_i$.  If $P_i$ is a lower saddle, then $P_i \cap h^{-1}(c - \eps)$ is a simple closed curve that bounds a disk $D_i$ in $S$ with exactly one minimum.  If $P_i$ is an upper saddle, there exists a similar disk $D_i$ with exactly one maximum.  In either case, $A_i = P_i \cup D_i$ is an annulus whose foliation contains exactly one saddle and one minimum or maximum, called a \emph{lower annulus} or an \emph{upper annulus}, respectively. \\

Now $\overline{T_k \setminus \cup A_i}$ is a collection of vertical annuli whose foliations contain no critical points, and we have decomposed $S$ into a collection of lower, upper, and vertical annuli, which we denote $\{A_1,\dots,A_r\}$.  In addition, these $r$ annuli as a collection have $r$ distinct boundary components, which we denote $\{\delta_1,\dots,\delta_r\}$.  Let $L$ be any level 2-sphere containing some $\delta_i$.  Of all curves in $T_k \cap L$ that are essential in $T_k$, consider a curve $\A$ that is innermost in $L$.  Thus $\A$ bounds a disk $D$ in $L$ containing no other simple closed curves in $T_k \cap L$.  Potentially, $D$ contains some inessential curves in $T_k \cap L$, but these curves bound disks in $T_k$, so after some gluing operations, we see that $\A$ bounds a compressing disk for $T_k$.  Since $T_k$ is compressible only on one side, it follows that $\A$ bounds a meridian disk of the solid torus $V_k$; hence, $\delta_i$ also bounds a meridian disk since it is parallel to $\A$ in $T_k$.  As $K$ is a $(p,q)$-cable, each $\delta_i$ has algebraic intersection $q$ with $k$, which implies $k$ has at least $q \cdot r$ intersections with the collection $\{\delta_1,\dots,\delta_r\}$. \\

We would like $k$ to be contained in the companion torus as efficiently as possible; for this purpose, we define efficient position, where we decompose $T_k$ into annuli as describe above.
\begin{defn}\label{effpos}
We say that $(k,T_k) \in \KT$ is an \textbf{efficient position} if $k \cap A_i$ is a collection of essential arcs in $A_i$ for every $i$.  Further, if $A_i$ is a vertical annulus, we require that each arc of $k \cap A_i$ contain no critical points, and if $A_i$ is an upper or lower annulus, we require that each arc of $k \cap A_i$ contains exactly one minimum or maximum.
\end{defn}
Note that if $(k,T_k)$ is an efficition position, it is implicit in Definition \ref{effpos} that the foliation of $T_k$ contains only essential saddles.  Of course, given an arbitrary element $(k,T_k) \in \KT$, we may not necessarily assume that $(k,T_k)$ is an efficient position, but we may employ the next lemma.
\begin{lemma}
\label{lem3}
For any $k \in \K$, there exists $(k',T_{k'}) \in \K$ such that $w(k') \leq w(k)$ and $(k',T_{k'})$ is an efficient position.
\begin{proof}

Let $(k,T_k) \in \KT$ such that $F_{T_k}$ contains no inessential saddles and decompose $T_k$ into annuli as described above.  Suppose that $\sum |k \cap \delta_i|$ is minimal up to isotopies that do not increase width.  By previous arguments, we suppose that after a series of type I and type II moves, we have for every vertical annulus $A_i$, $k \cap A_i$ consists of monotone essential arcs and inessential arcs with exactly one critical point, and for every lower or upper annulus $A_i$, $k \cap A_i$ consists of essential and inessential arcs with exactly one critical point.  As determined above, $\sum |k \cap \delta_i| \geq q \cdot r$.  If $\sum |k \cap \delta_i| = q \cdot r$, then every oriented intersection of $k$ with $\delta_i$ must occur with the same sign, so every arc of $k \cap A_i$ is essential and $(k,T_k)$ is an efficient position. \\

Conversely, if every arc of $k \cap A_i$ is essential for all $i$, then all oriented intersections of $k$ with $\delta_i$ occur with the same sign and $\sum |k \cap \delta_i| = q \cdot r$.  Thus, if $\sum |k \cap \delta_i| > q \cdot r$, there exists $A_j$ such that $k \cap A_j$ contains an inessential arc.  Let $\gamma$ be an inessential arc that is outermost in $A_j$, so that $\gamma$ cobounds a disk $\Delta$ with a level arc in $\pd A_j$ such that $\text{int}(\Delta) \cap k = \emp$.  Possibly after a disk slide if $A_j$ is upper or lower, we can slide $\gamma$ along $\Delta$ to remove two points of some $k \cap \delta_i$.  This isotopy does not increase the width of $k$, contradicting the minimality assumption above.

\end{proof}
\end{lemma}

\begin{figure}[h!]
  \centering
    \includegraphics[width=1.0\textwidth]{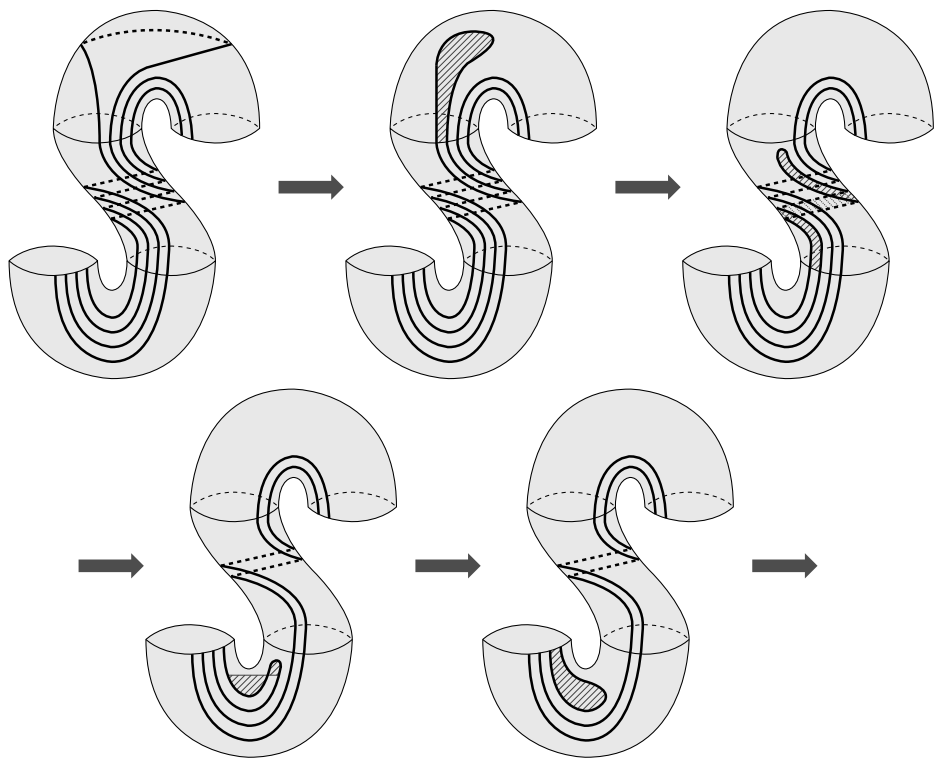}
    \caption{An example of the simplification suggested by Lemma \ref{lem3}.  In the first picture, an inessential arc $\gamma$ is contained in an upper annulus.  After a disk slide, we remove $\gamma$ from the upper annulus, and then remove it from the vertical annulus.  Now the inessential arc contains three critical points, so after a Type I move, we can remove it from the lower annulus, continuing as necessary until an efficient position is attained.}
\end{figure}

As a result of the proof of Lemma \ref{lem3}, if $(k,T_k)$ is an efficient position, we have that $k \cap A_i$ consists of $q$ essential arcs for each $i$.  See Figure 4.

\section{The width of the companion torus and the width of cable knots}
Neglecting $K$ for the moment, let $L$ be any knot in $S^3$ and $\LL$ the set of embeddings isotopic to $L$.  In this section, we define the width of a torus in $S^3$.  However, instead of modifying the standard definition that counts the number of intersections of with level 2-spheres, we will keep track of the order in which the upper and lower saddles occur in the foliation of the torus by $h$.  In the case of knots, observe that we can calculate the width of any embedding $l \in \LL$ if we know the order in which all minima and maxima of $l$ occur with respect to $h$.  To formalize this notion, let $Z$ be the free monoid generated by two elements, $m$ and $M$.  Now, define $\hat{Z}$ to be those elements $\sigma = m^{\A_1}M^{\n_1} \dots m^{\A_n}M^{\n_n}\in Z$ such that
\be
\item $\A_i,\n_i \neq 0$ for all $i$,
\item $\sum_{i=1}^j \A_i > \sum_{i=1}^j \n_i$ for all $j < n$, and
\item $\sum_{i=1}^n \A_i = \sum_{i=1}^n \n_i$.
\ee

We define a map from $\LL$ to $\hat{Z}$, $l \mapsto \sigma_l$, by the following:  Let $c_0 < \dots < c_p$ be the critical values of $h \mid_l$.  We create a word $\sigma_l$ consisting of $p+1$ letters by mapping the tuple $(c_0,\dots,c_p)$ a word by assigning $m$ to each minimum and $M$ to each maximum.  Next, we define the width of an element $\sigma = m^{\A_1}M^{\n_1} \dots m^{\A_n}M^{\n_n} \in \hat{Z}$ by
\[ w(\sigma) = 2\left( \sum \A_i \right)^2 - 4 \sum_{i > j} \A_i \n_j.\]
The following lemma should be expected:
\begin{lemma}\label{weq}
For any $l \in \LL$, $w(l) = w(\sigma_l)$.
\begin{proof}
Let $\sigma_l = m^{\A_1}M^{\n_1} \dots m^{\A_n}M^{\n_n}$.  Then the thick/thin tuple for $l$ is $(2\A_1,$ $2(\A_1-\n_1),2(\A_1-\n_1+\A_2),\dots,2(\A_1-\n_1 + \A_2 - \dots + \A_n))$.  From the width formula given by (\ref{wsquare}),
\begin{eqnarray*}
w(k) &=& \frac{1}{2}[ (2\A_1)^2 - (2(\A_1 - \n_1))^2 + (2(\A_1 - \n_1 + \A_2))^2 - \dots \\
&& \quad \dots +  (2(\A_1 - \n_1 + \A_2 - \dots + \A_n))^2 ] \\
&=& 2[\A_1^2 + \A_2^2 + 2\A_2(\A_1 - \n_1) + \dots \\
&& \quad \dots + \A_n^2 + 2\A_n(\A_1-\n_1+ \A_2 - \dots - \n_{n-1}) ] \\
&=& 2\left( \sum \A_i \right)^2 - 4 \sum_{i > j} \A_i \n_j,
\end{eqnarray*}
as desired.
\end{proof}
\end{lemma}
As an example, consider the embedding $l$ pictured in Figure 5.  A simple verification shows that $\sigma_l = m^3MmM^3$ and $w(l) = w(\sigma_l) = 28$. \\
\begin{figure}[h!]
  \centering
    \includegraphics[width=0.3\textwidth]{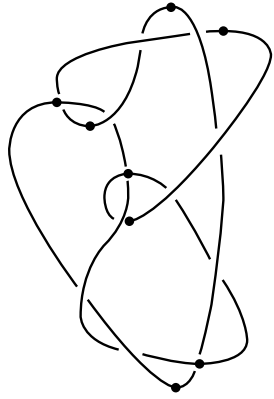}
    \caption{An embedding $l$ with $\sigma_l = m^3MmM^3$.}
\end{figure}

Next, we consider the collection of possible companion tori for a satellite knot with companion $L$.  Define $\Ss_L$ to be the collection of embedded tori $S$ such that
\be
\item $h\mid_S$ is Morse,
\item the foliation $F_S$ induced by $h$ contains no inessential saddles, and
\item $S$ is isotopic to the boundary of a regular neighborhood of $L$.
\ee
Now, as with $\LL$, we can define a map from $\Ss_L$ to $\hat{Z}$, $S \mapsto \sigma_S$, by the following:  Let $c_0 < \dots < c_p$ be the critical values corresponding to the saddles of $F_S$.  We create a word $\sigma_S$ consisting of $p+1$ letters by mapping the tuple $(c_0,\dots,c_p)$ to a word by assigning $m$ to each lower saddle and $M$ to each upper saddle.  It is clear that $\sigma_S \in Z$, and slightly more difficult to see that $\sigma_S \in \hat{Z}$.  The proof of this fact is left to the reader. \\

Using this association, we define the width of any torus $S \in \Ss_L$ by $w(S) = w(\sigma_S)$, and similar to the definition of knot width, we have
\begin{defn}
The \textbf{neighborhood width} of the knot $L$, $nw(L)$, is given by
\[ nw(L) = \min_{S \in \Ss_L} w(S).\]
\end{defn}
As a knot invariant, neighborhood width is not new; it is equivalent to knot width by the next lemma.
\begin{lemma}\label{lem4}
For any knot $L$, $nw(L) = w(L)$.
\begin{proof}
First, let $l_0$ be an embedding of $L$ such that $w(l_0) = w(L)$, and let $S_0$ be the boundary of a regular neighborhood of $l_0$ in $S^3$ such that every saddle of $S_0$ occurs slightly above a minimum or slightly below a maximum of $l_0$.  All such saddles are easily seen to be essential, so $S_0 \in \Ss_L$.  Further, $w(S_0) = w(\sigma_{S_0}) = w(\sigma_{l_0}) = w(L)$; hence
\[ nw(L) \leq w(L).\]

For the reverse inequality, let $S \in \Ss_L$.  Since $S$ is isotopic to a regular neighborhood of $L$, we have that any longitude of $S$ is isotopic to $L$.  Decompose $S$ into a collection of upper, lower, and vertical annuli $\{A_1,\dots,A_r\}$ as in Section \ref{eff}.  For each lower or upper annulus $A_i$, let $l_i$ be an essential arc in $A_i$ passing through the saddle with exactly one minimum or maximum.  For each vertical annulus $A_i$, let $l_i$ be a monotone arc connecting the two endpoints of the arcs contained in the annuli adjacent to $A_i$ (these must contain saddles).  Lastly, let $l$ be the union of all the arcs $l_i$, so that $l$ is a simple closed curve. \\

As shown above, the collection of boundary components $\{\delta_i\}$ of the annuli $\{A_i\}$ bound meridian disks for the solid torus bounded by $S$, and since $l$ intersects each $\delta_i$ once, $l$ is a longitude of $S$.  By construction $w(S) = w(\sigma_S) = w(\sigma_l) = w(l) \geq w(L)$, and since this is true for all $S \in \Ss_L$, we have
\[ nw(L) = w(L),\]
completing the proof. \\

\end{proof}
\end{lemma}
Finally, we have all the necessary tools to find the width of a cable knot.
\begin{thm}
Suppose that $K$ is an $(p,q)$-cable of a nontrivial knot $J$.  Then
\[ w(K) = q^2 \cdot w(J).\]
\begin{proof}
Let $j$ be an embedding of $J$ such that $w(j) = w(J)$, and let $T_j$ be the boundary of a regular neighborhood of $J$ such that every saddle of $T_j$ occurs slightly above a minimum or slightly below a maximum of $j$.  Suppose that $\sigma_j = m^{\A_1}M^{\n_1} \dots$ $m^{\A_n}M^{\n_n}$.  We can cable $j$ along $T_j$, creating an embedding $k \in \K$ such that $\sigma_k = m^{q\A_1}M^{q\n_1} \dots m^{q\A_n}M^{q\n_n}$.  By Lemma \ref{weq},
\[ w(K) \leq  w(\sigma_k) = q^2 \cdot w(\sigma_j) = q^2 \cdot w(J).\]
We call $k$ an ``obvious" cabling of $j$. \\

On the other hand, suppose $k' \in \K$ is a given thin position.  By Lemmas \ref{lem2} and \ref{lem3}, there exists a torus $T_{k'}$ such that $(k',T_{k'}) \in \KT$ is an efficient position.  Let $\sigma_{T_{k'}} = m^{\A'_1}M^{\n'_1} \dots m^{\A'_{n'}}M^{\n'_{n'}}$.  Note that $T_{k'} \in \Ss_J$, and thus $w(T_{k'}) \geq w(J)$ by Lemma \ref{lem4}.  We decompose $T_{k'}$ into a collection of upper, lower, and vertical annuli $\{A_1,\dots,A_r\}$ as above in Section \ref{eff}.  For each upper annulus $A_i$, $k \cap A_i$ consists of $q$ essential arcs, each containing one maximum.  Suppose $c_i$ is the critical value corresponding to the saddle in $A_i$.  Then there is an isotopy of $k$ supported in $A_i$ taking the $q$ arcs to arcs in $h^{-1}((c_i,c_i+\eps])$, and this isotopy does not increase $w(k)$.   A similar statement is true for each lower annulus $A_i$.  As all critical points of $k$ are contained in upper or lower annuli, we can compute $\sigma_k = m^{q\A'_1}M^{q\n'_1} \dots m^{q\A'_{n'}}M^{q\n'_{n'}}$, and thus
\[ w(k) = w(m^{q\A'_1}M^{q\n'_1} \dots m^{q\A'_{n'}}M^{q\n'_{n'}}) = q^2 \cdot w(T_{k'}) \geq q^2 \cdot w(J),\]
completing the proof.

\end{proof}
\end{thm}
The proof of the theorem reveals that $k$ is a thin position of $K$ if and only if it is an ``obvious" cable of a thin position of $J$.
\begin{cor}
For a knot $J$ in $S^3$, the following are equivalent:
\be
\item[(a)] $J$ is bridge-thin.
\item[(b)] Every cable of $J$ is bridge-thin.
\item[(c)] There exists a cable of $J$ that is bridge-thin.
\ee
\end{cor}
 
\section{Meridional smallness of cable knots}\label{small}
In the previous section, we showed that bridge-thinness is a property preserved by cabling.  Here we wish to show the same is true for meridional smallness and mp-smallness, defined below.  We need several other definitions first.
\begin{defn}
Suppose $S \subset M$ is a properly embedded surface.  A \textbf{compressing disk} for $S$ in $M$ is an embedded disk $D$ such that $\pd D \subset S$, $\text{int}(D) \cap S = \emp$, and $\pd D$ does not bound a disk in $S$.  We say $S$ is \textbf{incompressible} if there does not exist a compressing disk for $S$ in $M$.
\end{defn}
\begin{defn}
Suppose $S \subset M$ is a properly embedded surface, with $\pd S \neq \emp$ and $\pd S \subset \pd M$.  A \textbf{$\pd$-compressing disk} $D$ for $S$ in $M$ is an embedded disk $D$ such that $\pd D$ consists of two arcs $\A$ and $\n$, where $\A \subset \pd M$ and $\n \subset S$, $\text{int}(D) \cap S = \emp$, and such that $\n$ does not cobound a disk in $S$ with an arc in $\pd S$.  We say $S$ is \textbf{$\pd$-incompressible} if there does not exist a $\pd$-compressing disk for $S$ in $M$.
\end{defn}
\begin{defn}
A surface $S \subset M$ is \textbf{essential} if it is incompressible, $\pd$- incompressible, and not parallel to $\pd M$.
\end{defn}
Finally, we arrive at
\begin{defn}
A knot $L \subset S^3$ is \textbf{meridionally small} if $E(L)$ contains no essential surface $S$ with $\pd S$ consisting of meridian curves of $\overline{\eta(K)}$, where $\eta$ denotes a regular neighborhood.  A knot $L \subset S^3$ is \textbf{meridionally planar small}, or mp-small, if $E(L)$ contains no essential planar surface $S$ with meridional boundary.
\end{defn}
Note that if $S$ is a properly embedded surface with meridional boundary in a knot exterior $E(L)$ and $S$ is $\pd$-compressible, then either $S$ is compressible or $\pd$-parallel.  In addition, all components of $\pd E(L) \setminus \eta(\pd S)$ are annuli, so if $S$ is $\pd$-parallel, $S$ must be an annulus.   Thus to show $S$ is essential it suffices to show $S$ is incompressible and not a $\pd$-parallel annulus. \\ 

Recall that $K$ is a $(p,q)$-cable of $J$, with pattern $\hat{K}$ a $(p,q)$-torus knot contained in a solid torus $V$.  Letting $C_{p,q} = V \setminus \eta(\hat{K})$, we observe that we can decompose $E(K)$ as $E(J) \cup C_{p,q}$, where the attaching map depends on the framing $\varphi:V \rightarrow S^3$ that maps a core of $V$ to $J$.  Following \cite{lither}, we call $C_{p,q}$ a \emph{$(p,q)$-cable space}.  See Figure 6.
\begin{figure}[h!]
  \centering
    \includegraphics[width=0.6\textwidth]{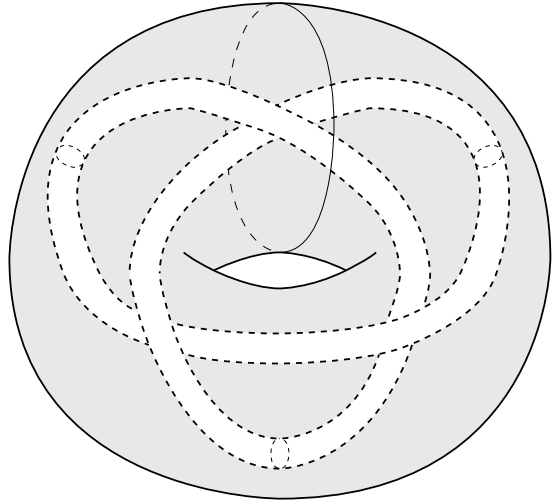}
    \caption{The cable space $C_{3,2}$.}
\end{figure}
Note that $C_{p,q}$ is Seifert fibered with one exceptional fiber, and it has two torus boundary components, $\pd V$ (the \emph{outer boundary}, denoted $\pd_+ C_{p,q}$) and $\pd \, \eta(\hat{K})$ (the \emph{inner boundary}, denoted $\pd_-C_{p,q}$).  In order to understand essential meridional surfaces in $E(K)$, we require the next lemma.  We note that this lemma is implied by Lemma 3.1 of \cite{lither}, but we do not need its full generality and so we provide an elementary proof.  For a treatment of Seifert fibered spaces, see \cite{hatcher}.  A surface in a Seifert fibered space is \emph{horizontal} if it is transverse to all fibers and \emph{vertical} if it is a union of fibers.
\begin{lemma}\label{lem5}
Suppose $S$ is a connected incompressible surface in $C_{p,q}$, where $S \cap \pd_- C_{p,q} \neq \emp$ and each component of $S \cap \pd_- C_{p,q}$ is a meridian curve.  Then $S$ is a disk with $q$ punctures, and $S \cap \pd_+ C_{p,q}$ is a meridian curve of $V$.
\begin{proof}
Since $C_{p,q}$ is Seifert fibered, every incompressible surface is either horizontal or vertical by \cite{hatcher}.  By  assumption $S$ has meridional boundary in $\pd_- C_{p,q}$, so $S$ is horizontal.  We can extend $S$ to a horizontal surface $S'$ in a Seifert fibered solid torus by attaching meridian disks to each component of $S \cap \pd_- C_{p,q}$ and gluing $\eta(\hat{K})$ back into $C_{p,q}$.  But every connected horizontal surface in such a torus is a meridian disk, so $S'$ is a meridian disk, from which it follows that $S$ is a disk with $q$ punctures and $\pd S \cap \pd_+ C_{p,q} = \pd S'$ is a meridian curve.
\end{proof}
\end{lemma}
On the other hand, let $S'$ be a meridian disk of $V$ intersecting $\hat{K}$ minimally, with $S = S' \cap C_{p,q}$.  Then $C_{p,q}$ cut along $S$ is homeomorphic to $S \X I$, from which it follows that $S$ is incompressible in $C_{p,q}$.  We use these facts in the proof of the next theorem.
\begin{thm}
For a knot $J$ in $S^3$, the following are equivalent:
\be
\item[(a)] $J$ is meridionally small.
\item[(b)] Every cable of $J$ is meridionally small.
\item[(c)] There exists a cable of $J$ that is meridionally small.
\ee
\begin{proof}
Suppose that $E(J)$ contains an essential meridional surface $S$, and let $K$ be a $(p,q)$-cable of $J$, so that $E(K) = E(J) \cup_T C_{p,q}$, where $T = \pd E(J)$.  Then every component of $\pd S$ bounds a disk with $q$ punctures in $C_{p,q}$, and thus we can construct a meridional surface $S'$ in $E(K)$ by gluing such a disk $D$ to each boundary component of $\pd S$.  We claim that $S'$ is essential in $E(K)$.  By the above, $S' \cap C_{p,q}$ is incompressible in $C_{p,q}$.  Note that $|\pd S'| = q \cdot |\pd S| > 2$, so $S'$ is not a $\pd$-parallel annulus.  Suppose $\Delta$ is a compressing disk for $S'$ in $E(K)$.  If $\Delta \cap T = \emp$, then either $\Delta \subset C_{p,q}$, which is ruled out by the argument above, or $\Delta \subset E(J)$, which contradicts the incompressibility of $S$. \\

Thus, $\Delta \cap T \neq \emp$, and suppose $\Delta$ is chosen so that $|\Delta \cap T|$ is minimal.  Since $T$ is incompressible in $E(K)$, $\Delta \cap T$ contains no simple closed curves.  Let $\A$ be an arc in $\Delta \cap T$ that is outermost in $\Delta$, so that $\A$ cobounds a disk $\Delta' \subset \Delta$ with an arc $\n \subset \pd \Delta$ such that $\text{int}(\Delta') \cap T = \emp$.  There are two cases to consider:  First, suppose that $\Delta' \subset C_{p,q}$.  Then both endpoints of $\A$ are contained in the same component of $S' \cap T$, which means that $\A$ is inessential in an annular component of $T$ and cobounds a disk $\Delta'' \subset T$ with an arc $\gamma \subset S' \cap T$.  Gluing $\Delta'$ to $\Delta''$ along $\A$ yields a disk $\Delta^*$ with $\pd \Delta^* \subset S' \cap C_{p,q}$, and by the incompressibility of $S' \cap C_{p,q}$ in $C_{p,q}$, $\pd \Delta^* = \n \cup \gamma$ bounds a disk $D$ contained in $S' \cap C_{p,q}$.  Sliding $\n$ along $D$, we can remove at least one intersection of $\Delta$ with $T$, contradicting the minimality of $|\Delta \cap T|$. \\

In the second case, suppose that $\Delta' \subset E(K)$.  If $\A$ is an inessential arc in some annular component of $T$, the above argument holds.  If $\A$ is essential, then $\Delta'$ is a $\pd$-compressing disk for $S$ in $E(J)$, contradicting the assumption that $S$ is essential.  We conclude that $S'$ is an essential meridional surface in $E(K)$, showing that (c) implies (a). \\

Now, suppose that there exists a $(p,q)$-cable $K$ of $J$ such that $E(K)$ contains an essential meridional surface $R$.  As above, $E(K) = E(J) \cup_T C_{p,q}$.  Assume $|R \cap T|$ is minimal.  We claim that $R \cap C_{p,q}$ is incompressible in $C_{p,q}$.  Suppose not.  Then there is a compressing disk $\Delta \subset C_{p,q}$ for $R \cap C_{p,q}$.  Since $R$ is incompressible in $E(K)$, $\pd \Delta$ bounds a disk $\Delta' \subset R$, and by isotopy we can replace $\Delta'$ with $\Delta$ and reduce the number of intersections of $R$ with $T$, a contradiction.  By Lemma \ref{lem5}, each component of $R \cap C_{p,q}$ is a punctured disk, and it follows that $R' = R \cap E(J)$ is a meridional surface.  If $R'$ is a $\pd$-parallel annulus in $E(J)$, we can find a $\pd$-compressing disk for $R'$ in $E(J)$, and any compressing or $\pd$-compressing disk for $R'$ in $E(J)$ yields a compressing or $\pd$-compressing disk for $R$ in $E(K)$, a contradiction.  Thus, $R'$ is an essential meridional surface in $E(J)$, and (a) implies (b).  Clearly, (b) implies (c), completing the proof.

\end{proof}
\end{thm}
It should be noted that in the above proof, if either surface $S$ or $R$ is planar, then the induced surface $S'$ or $R'$ is also planar, yielding:
\begin{thm}\label{mpsmall}
For a knot $J$ in $S^3$, the following are equivalent:
\be
\item[(a)] $J$ is mp-small.
\item[(b)] Every cable of $J$ is mp-small.
\item[(c)] There exists a cable of $J$ that is mp-small.
\ee
\end{thm}

\section{Destabilizing non-minimal bridge positions of cables}
If $l$ is a bridge position for a knot $L$, then $l$ has exactly one thick level $\hat{A} = h^{-1}(a_1)$.  We call this level a bridge sphere and say that two bridge positions $l$ and $l'$ equivalent if there exists an isotopy $f_t$, called a \emph{bridge isotopy}, taking $l$ to $l'$ such that for all $t$, $f_t(\hat{A})$ is a bridge sphere for $f_t(l)$.  If $l$ is a bridge position admitting a type I move, we say that $l$ is \emph{stabilized}.  In this context a type I move is also called a \emph{destabilization}.  Since a type I move cancels a minimum and a maximum, it is clear that every minimal bridge position is not stabilized.  For $K$ a $(p,q)$-cable of $J$, we have the following:
\begin{lemma}\label{stab1}
Suppose that $k$ is a bridge position.  Then there exists $(k',T_{k'}) \in \KT$ such that $k'$ is equivalent to $k$ and $(k',T_{k'})$ is an efficient position, or $k$ is stabilized.
\begin{proof}
Let $(k,T_k) \in \KT$.  By Lemma 2.4 of \cite{ozawa}, there exists $(k',T_{k'}) \in \KT$ such that $k'$ is equivalent to $k$ and all saddles in the foliation of $T_{k'}$ by $h$ are essential.  Note that any disk slide of $k'$ is supported in a neighborhood away from a bridge sphere; thus, if $k''$ is the result of a disk slide on $k'$, $k''$ is equivalent to $k'$.  Now, by Lemma \ref{lem3}, either $(k',T_{k'})$ is an efficient position, or $k'$ admits a destabilization, possibly after a disk slide.
\end{proof}
\end{lemma}

In \cite{ozawa}, Makoto Ozawa shows that if $K$ is a torus knot, then every non-minimal bridge position of $K$ is stabilized.  We essentially follow his proof, which we will summarize here, to show that if $J$ is mp-small and every non-minimal bridge position of $J$ is stabilized, then every cable $K$ of $J$ has the same property.  Ozawa utilizes the following theorem, proved by Hayashi and Shimokawa \cite{hayashi} and later by Tomova \cite{tomova}:
\begin{thm}
If $k \in \K$ is a bridge position admitting a type II move, then either $k$ is stabilized or $E(K)$ contains an essential meridional planar surface.
\end{thm}
Thus, if $K$ is mp-small, any bridge position admitting a type II move is stabilized. \\

From this point on, suppose $J$ is mp-small.  By Theorem \ref{mpsmall}, $K$ is also mp-small.  We recall that for $(k,T_k) \in \KT$, the solid torus bounded by $T_k$ is denoted $V_k$.  Although \cite{ozawa} deals specifically with torus knots, some results apply in the setting of mp-small cable knots.  Lemmas 3.2, 3.3, and 3.4 of \cite{ozawa} together imply
\begin{lemma}\label{stab3}
Suppose $k$ is a bridge position.  Then there exists $(k',T_{k'}) \in \KT$ an efficient position such that $k'$ is equivalent to $k$ and a bridge sphere $\hat{A} = h^{-1}(a_1)$ such that all upper saddles of $T_{k'}$ occur above $\hat{A}$, all lower saddles occur below $\hat{A}$, and $V_{k'} \cap \hat{A}$ is a collection of disks, or $k$ is stabilized.
\end{lemma}
Thus, $h$ foliates the solid torus $V_{k'}$ by disks, and any such $(k',T_{k'})$ induces a bridge position $j$ of $J$ by taking $j$  a longitude of $T_{k'}$ as constructed in Lemma \ref{weq}.  This fact also ensures that any upper or lower disk for $j$ can be chosen to miss $\text{int}(V_{k'})$ and thus can easily be modified to an upper or lower disk for $k'$.  Finally, we have
\begin{thm}
Suppose $K$ is a $(p,q)$-cable with companion $J$, where $J$ is mp-small.
\be
\item[(a)] If every non-minimal bridge position of $J$ is stabilized, then every non-minimal bridge position of $K$ is stabilized.
\item[(b)] The cardinality of the collection of minimal bridge positions of $K$ does not exceed the cardinality of the collection of minimal bridge positions of $J$.
\ee
\begin{proof}
(a) Suppose $k \in \K$ is a bridge position.  By Lemmas \ref{stab1} and \ref{stab3}, we  may pass to an equivalent bridge position and assume there exists $(k,T_{k}) \in \KT$ an efficient position such that $h$ foliates $V_k$ by disks, or else $k$ is stabilized.  Let $\hat{A}$ be a bridge sphere for $k$ guaranteed by Lemma \ref{stab3}, and note that the induced bridge position $j$ of $J$ satisfies $b(j) = \frac{1}{2} |\hat{A} \cap T_k|$.  There are two cases to consider:  First, suppose that $|\hat{A} \cap T_k| = 2 \cdot b(J)$.  Since each upper annulus contains $q$ maxima, we have $b(k) = q \cdot b(J)$.  By Theorem 1 of \cite{bridge}, $b(K) \geq q \cdot b(J)$, and thus $k$ is a minimal bridge position of $K$.  On the other hand, suppose that $|\hat{A} \cap T_k| > 2 \cdot b(J)$.  Then $j$ is a non-minimal bridge position for $J$, and by assumption $j$ is stabilized.  It follows that $k$ is stabilized, finishing the first part of the proof. \\

(b) For each minimal bridge position $j$ of $J$, assign an obvious cabling $k \subset \pd \eta(j)$ of $j$ with $q \cdot b(J)$ maxima.  It suffices to prove that the association $\varphi: j \mapsto k$ is surjective.  First, suppose that $k'$ is another obvious cabling of a bridge position of $j'$ equivalent to $j$ with $q \cdot b(J)$ maxima.  Then there is a bridge isotopy taking $j'$ to $j$, which induces an isotopy from $\pd \eta(j')$ to $\pd \eta(j)$ which is also a bridge isotopy of $k'$.  Since $k$ is isotopic to $k'$ in $\pd \eta(j)$, we can remove intersections of $k$ and $k'$ and then push $k'$ onto $k$ by a bridge isotopy, and so $\varphi$ is well-defined.  Now, if $k$ is a minimal bridge position of $K$, we may pass to an equivalent bridge position and assume there exists $(k,T_k)$ such that $h$ foliates $V_k$ by disks by Lemma \ref{stab3}.  Thus, the induced minimal bridge position $j$ for $J$ from Lemma \ref{weq} satisfies $\varphi(j) = k$, as desired.

\end{proof}
\end{thm}

From \cite{otal} and \cite{ozawa}, we have
\begin{cor}
If $K$ is an $n$-fold cable of a torus knot or a 2-bridge knot, then any non-minimal bridge position of $K$ is stabilized.  Additionally, if $K$ is an $n$-fold cable of a torus knot, it has a unique minimal bridge position, and if $K$ is an $n$-fold cable of a 2-bridge knot, it has at most two minimal bridge positions.
\end{cor}
In \cite{coward}, Coward proves that if $J$ is a hyperbolic knot, then $J$ has finitely many minimal bridge positions.  Hence
\begin{cor}
If $K$ is a $(p,q)$-cable of $J$, where $J$ is hyperbolic and mp-small, then $K$ has finitely many minimal bridge positions.
\end{cor}

\end{document}